\def\Bbb#1{{\fam\msbfam\relax#1}}
\font\fivemsb=msbm5
\font\sevenmsb=msbm7
\font\tenmsb=msbm10
\newfam\msbfam
\textfont\msbfam=\tenmsb
\scriptfont\msbfam\sevenmsb
\scriptscriptfont\msbfam\fivemsb

\def\spc{{\Bbb C}}
\def\vf{\varphi}
\def\wt{\widetilde}
\def\wh{\widehat}
\def\a{\alpha}
\def\b{\beta}

\def\F{\Phi}

\def\g{\gamma}
\def\m{\mu}
\def\r{\rho}

\def\p{\pi}

\def\x{\xi}
\def\y{\eta}

\def\ce{{\cal E}}

\def\la{\langle}
\def\ra{\rangle}
\def\lla{\langle\!\langle}
\def\rra{\rangle\!\rangle}

\documentstyle[verbatim]{amsart}

\theoremstyle{plain}
\newtheorem{theorem}{Theorem}[section]

\newtheorem{lemma}[theorem]{Lemma}

\theoremstyle{definition}
\newtheorem{example}[theorem]{Example}

\numberwithin{equation}{section}

\begin{document}
\title[Characterization of test functions]%
   {Characterization of test functions in CKS-space}

\author{Nobuhiro Asai}
\author{Izumi Kubo}
\author{Hui-Hsiung Kuo}
\address{Nobuhiro Asai: Graduate School of Mathematics\\
  Nagoya University \\ Nagoya 464-8602 \\ JAPAN}
\address{Izumi Kubo: Department of Mathematics \\
Faculty of Science \\ Hiroshima University \\
Higashi-Hiroshima 739-8526 \\ JAPAN}
\address{Hui-Hsiung Kuo: Department of Mathematics\\
  Louisiana State University \\ Baton Rouge \\
LA 70803 \\ USA}

\thanks{Research supported by the Daiko Foundation 1998
(N.A.) and Academic Frontier in Science of Meijo
University (H.-H.K.)}

\maketitle

\begin{abstract}
We prove a characterization theorem for the test functions
in a CKS-space. Some crucial ideas concerning the growth
condition are given.

\end{abstract}

\section{Introduction} \label{sec:1}

Let $\ce$ be a real countably-Hilbert space with topology
given by a sequence of norms $\{|\cdot|\}_{p=0}^{\infty}$
(see \cite{kuo96}.) Let $\ce_{p}$ be the completion of $\ce$
with respect to the norm $|\cdot|_{p}$. Assume the following
conditions:

\begin{itemize}
\item[(a)] There exists a constant $0<\r<1$ such that
$|\cdot|_{0} \leq \r |\cdot|_{1} \leq \cdots \leq
\r^{p} |\cdot|_{p} \leq \cdots$.
\item[(b)] For any $p\geq 0$, there exists $q>p$ such that
the inclusion map $i_{q, p}: \ce_{q} \to \ce_{p}$ is a
Hilbert-Schmidt operator and $\|i_{q, p}\|_{HS} < 1.$
\end{itemize}

By using the Riesz representation theorem to identify
$\ce_{0}$ with its dual space we get the continuous
inclusion maps:
\begin{equation}
\ce \subset \ce_{p} \subset \ce_{0} \subset \ce_{p}'
\subset \ce', \qquad p\geq 0,  \notag
\end{equation}
where $\ce_{p}'$ and $\ce'$ are the dual spaces of $\ce_{p}$
and $\ce$, respectively. Condition (b) says that $\ce$ is
a nuclear space and so $\ce \subset \ce_{0} \subset \ce'$
is a Gel'fand triple. Let $\m$ be the standard Gaussian
measure on $\ce'$. Let $(L^{2})$ denote the Hilbert space of
$\m$-square integrable functions on $\ce'$. By the
Wiener-It\^o theorem each $\vf$ in $(L^{2})$ can be uniquely
expressed as
\begin{equation} \label{eq:1-1}
\vf = \sum_{n=0}^{\infty} I_{n}(f_{n}) = \sum_{n=0}^{\infty}
\la :\!\cdot^{\otimes n}\!:, f_{n}\ra, \qquad
f_{n}\in \ce_{0}^{\wh\otimes n},
\end{equation}
and the $(L^{2})$-norm $\|\vf\|_{0}$ of $\vf$ is given by
\begin{equation}
\|\vf\|_{0} = \left(\sum_{n=0}^{\infty} n!|f_{n}|_{0}^{2}
\right)^{1/2}.   \notag
\end{equation}

Now, we describe the spaces of test and generalized functions
on the space $\ce'$ introduced by Cochran et al.~in a recent
paper \cite{cks}. Let $\{\a(n)\}_{n=0}^{\infty}$ be a
sequence of numbers satisfying the following conditions:
\begin{itemize}
\item[(A1)] $\a(0)=1$ and $\inf_{n\geq 0}\,\a(n) > 0$.
\item[(A2)] $\lim_{n\to\infty} \left({\a(n)
\over n!}\right)^{1/n} =0$.
\end{itemize}

Let $\vf\in (L^{2})$ be represented as in Equation
(\ref{eq:1-1}). For each nonnegative integer $p$, define
\begin{equation} \label{eq:a}
\|\vf\|_{p, \a} = \left(\sum_{n=0}^{\infty} n! \a(n)
|f_{n}|_{p}^{2} \right)^{1/2}.
\end{equation}
Let $[\ce_{p}]_{\a} = \{\vf\in (L^{2})\, ;\, \|\vf\|_{p, \a}
<\infty\}$. Define the space $[\ce]_{\a}$ of test functions
to be the projective limit of the family $\{[\ce_{p}]_{\a};\,
p\geq 0\}$. Its dual space $[\ce]_{\a}^{*}$ is the space of
generalized functions. By identifying $(L^{2})$ with its
dual we get the following continuous inclusion maps:
\begin{equation}
[\ce]_{\a} \subset [\ce_{p}]_{\a} \subset (L^{2}) \subset
[\ce_{p}]_{\a}^{*} \subset [\ce]_{\a}^{*},
\quad p\geq 0. \notag
\end{equation}
If $\F\in [\ce_{p}]_{\a}^{*}$ (the dual space of
$[\ce_{p}]_{\a}$) is represented by $\F =
\sum_{n=0}^{\infty} \la :\!\cdot^{\otimes n}\!:, F_{n}\ra$,
then its $[\ce_{p}]_{\a}^{*}$-norm is given by
\begin{equation}  \label{eq:a-1}
\|\vf\|_{-p, 1/\a} = \left(\sum_{n=0}^{\infty} {n! \over
\a(n)} |F_{n}|_{-p}^{2} \right)^{1/2}.
\end{equation}

This Gel'fand triple $[\ce]_{\a} \subset (L^{2}) \subset
[\ce]_{\a}^{*}$ is called the CKS-space associated with a
sequence $\{\a(n)\}_{n=0}^{\infty}$ of numbers satisfying
the above conditions (A1) and (A2).

Several characterization theorems for generalized functions
in $[\ce]_{\a}^{*}$ have been proved in the paper \cite{cks}.
However, no characterization theorem for test functions in
$[\ce]_{\a}$ is given. The purpose of the present paper is
to prove such a theorem. In addition we will mention some
crucial ideas in order to get a complete description of
the characterization theorems for test and generalized
functions in our ongoing research collaboration project.
We remark that similar results have been obtained by
Gannoun et al.~\cite{gho}.

\section{Characterization theorems} \label{sec:2}

For $\x\in \ce_{c}$ (the complexification of $\ce$,) the
renormalized exponential function $:\!e^{\la \cdot,
\x\ra}\!:$ is defined by
\begin{equation}  \label{eq:a-2}
:\!e^{\la \cdot, \x\ra}\!: \,= \sum_{n=0}^{\infty}
{1\over n!}\, \la :\!\cdot^{\otimes n}\!:,
\x^{\otimes n}\ra.
\end{equation}
For any $p\geq 0$, we have
\begin{equation} \label{eq:2-1}
\|\!:\!e^{\la \cdot, x\ra}\!:\!\|_{p, \a} =
   G_{\a}\big(|\x|_{p}^{2}\big)^{1/2},
\end{equation}
where $G_{\a}$ is the exponential generating function of
the sequence $\{\a(n)\}_{n=0}^{\infty}$, i.e.,
\begin{equation}
G_{\a}(z) = \sum_{n=0}^{\infty} {\a(n) \over n!} z^{n},
\qquad z\in\spc.  \notag
\end{equation}
By condition (A2) of the sequence $\{\a(n)\}_{n=0}^{\infty}$,
the function $G_{\a}$ is an entire function. Equation
(\ref{eq:2-1}) implies that $:\!e^{\la \cdot, x\ra}\!:$ is
a test function in $[\ce]_{\a}$ for any $\x\in\ce_{c}$.

The $S$-transform of a generalized function $\F$ in
$[\ce]_{\a}^{*}$ is the function $S\F$ defined on
$\ce_{c}$ by
\begin{equation}
S\F(\x) = \lla \F, :\!e^{\la \cdot, x\ra}\!:\rra,
\qquad \x\in\ce_{c},  \notag
\end{equation}
where $\lla \cdot, \cdot\rra$ is the bilinear pairing of
$[\ce]_{\a}^{*}$ and $[\ce]_{\a}$.

We state three conditions on the sequence
$\{\a(n)\}_{n=0}^{\infty}$ as follows:

\begin{itemize}
\item[(B1)] $\limsup_{n\to\infty}\left({n!\over \a(n)}
\inf_{r>0} {G_{\a}(r) \over r^{n}}\right)^{1/n} < \infty.$
\item[(B2)] The sequence $\g(n)={\a(n) \over n!}, n\geq 0$,
is log-concave, i.e.,
\begin{equation}
\g(n)\g(n+2) \leq \g(n+1)^{2}, \qquad\forall n\geq 0. \notag
\end{equation}
\item[(B3)] The sequence $\{\a(n)\}_{n=0}^{\infty}$ is
log-convex, i.e.,
\begin{equation}
\a(n)\a(n+2) \geq \a(n+1)^{2}, \qquad\forall n\geq 0. \notag
\end{equation}
\end{itemize}

Condition (B1) is used in the characterization theorem for
generalized functions. By Corollary 4.4 in \cite{cks}
condition (B2) implies condition (B1). Condition (B3)
implies condition ($\wt{\text{B}}2$) to be defined below.

The following characterization theorem for generalized
functions in $[\ce]_{\a}^{*}$ has been proved in \cite{cks}.

\begin{theorem} \label{thm:1}
Suppose the sequence $\{\a(n)\}_{n=0}^{\infty}$ satisfies
conditions {\em (A1)} and {\em (A2)}. {\em (I)} Let $\F\in
[\ce]_{\a}^{*}$. Then the $S$-transform $F=S\F$ of $\F$
satisfies the conditions:
\begin{itemize}
\item[(1)] For any $\x, \y$ in $\ce_{c}$, the function
$F(z\x+\y)$ is an entire function of $z\in\spc$.
\item[(2)] There exist constants $K, a, p \geq 0$ such that
\begin{equation} \label{eq:d}
|F(\x)| \leq K G_{\a} \big(a|\x|_{p}^{2}\big)^{1/2},
  \qquad \x\in\ce_{c}.
\end{equation}
\end{itemize}
{\em (II)} Conversely, assume that condition {\em (B1)}
holds and let $F$ be a function on $\ce_{c}$ satisfying
the above conditions {\em (1)} and {\em (2)}. Then there
exists a unique $\F\in [\ce]_{\a}^{*}$ such that $F=S\F$.
\end{theorem}

We mention that condition (2) is actually equivalent to the
condition: there exist constants $K, p \geq 0$ such that
\begin{equation}
|F(\x)| \leq K G_{\a} \big(|\x|_{p}^{2}\big)^{1/2},
  \qquad \x\in\ce_{c}.  \notag
\end{equation}
This fact can be easily checked by using the fact that
$|\x|_{p} \leq \r^{q-p} |\x|_{q}$ for any $q\geq p$.
Having the constant $a$ in condition (2) is for convenience
to check whether a given function $F$ satisfies the
condition.

By this theorem, if we assume condition (B1), then
conditions (1) and (2) are necessary and sufficient for a
function $F$ defined on $\ce_{c}$ to be the $S$-transform
of a generalized function in $[\ce]_{\a}^{*}$. As mentioned
above, condition (B2) implies condition (B1). Thus under
condition (B2), conditions (1) and (2) are also necessary
and sufficient for $F$ to be the $S$-transform of a
generalized function in $[\ce]_{\a}^{*}$.

For the characterization theorem for test functions in
$[\ce]_{\a}$, we need to define the exponential generating
function of the sequence $\{{1\over \a(n)}\}_{n=0}^{\infty}$:
\begin{equation} \label{eq:b}
G_{1/\a}(z) = \sum_{n=0}^{\infty} {1\over n!\a(n)} z^{n},
  \qquad z\in\spc,
\end{equation}
Moreover, we need the corresponding conditions (A2) (B1)
(B2) for the sequence $\{{1 \over \a(n)}\}_{n=0}^{\infty}$:

\begin{itemize}
\item[($\wt{\text{A}}2$)] $\lim_{n\to\infty} \left(
{1 \over n!\a(n)}\right)^{1/n} =0$.
\item[($\wt{\text{B}}1$)] $\limsup_{n\to\infty}\left(n!\a(n)
\inf_{r>0} {G_{1/\a}(r) \over r^{n}}\right)^{1/n} < \infty.$
\item[($\wt{\text{B}}2$)] The sequence $\{{1\over n!\a(n)}
\}_{n=0}^{\infty}$ is log-concave.
\end{itemize}

It follows from condition ($\wt{\text{A}}2$) that the
exponential generating function $G_{1/\a}$ is an entire
function. Note that by condition (A1), $\a(n)\geq \a_{0}$
for all $n$, where $\a_{0}=\inf_{n\geq 0} \a(n)$. Hence
by the Stirling formula
\begin{equation}
\left({1 \over n!\a(n)}\right)^{1/n}
\leq \left({1 \over n!\a_{0}}\right)^{1/n}
\sim \left({1\over \a_{0}\,\sqrt{2\p n}}\right)^{1/n}\,
  {e\over n} \to 0, \quad \text{as~} n\to\infty. \notag
\end{equation}

This shows that condition (A1) implies condition
($\wt{\text{A}}2$). On the other hand, by applying
Corollary 4.4 in \cite{cks} to the sequence
$\{{1\over \a(n)}\}_{n=0}^{\infty}$ we see that condition
($\wt{\text{B}}2$) implies condition ($\wt{\text{B}}1$).
Moreover, it is easy to check that condition (B3) implies
condition ($\wt{\text{B}}2$).

Now, we study the characterization theorem for test
functions in $[\ce]_{\a}$. First we prove a lemma.

\begin{lemma} \label{lem:1}
Assume that condition {\em ($\wt{\text{B}}1$)} holds and
let $F$ be a function on $\ce_{c}$ satisfying the
conditions:
\begin{itemize}
\item[(1)] For any $\x, \y$ in $\ce_{c}$, the function
$F(z\x+\y)$ is an entire function of $z\in\spc$.
\item[(2)] There exist constants $K, a, p \geq 0$ such that
\begin{equation}
|F(\x)| \leq K G_{1/\a} \big(a|\x|_{-p}^{2}\big)^{1/2},
  \qquad \x\in\ce_{c}.  \notag
\end{equation}
\end{itemize}
Then there exists a unique $\vf\in [\ce]_{\a}^{*}$ such
that $F=S\vf$. In fact, $\vf\in [\ce_{q}]_{\a}$ for any
$q\in [0, p]$ satisfying the condition
\begin{equation} \label{eq:2-2}
ae^{2}\|i_{p, q}\|_{HS}^{2}\limsup_{n\to\infty}\left(
n!\a(n) \inf_{r>0}{G_{1/\a}(r) \over r^{n}}\right)^{1/n}
< 1
\end{equation}
and
\begin{equation} \label{eq:2-3}
\|\vf\|_{q, \a}^{2} \leq {K^{2} \over 2\p} \sum_{n=0}^{\infty}
\left(n!\a(n)\inf_{r>0}{G_{1/\a}(r) \over r^{n}}\right)
\Big(ae^{2}\|i_{p, q}\|_{HS}^{2}\Big)^{n}.
\end{equation}
\end{lemma}

\begin{pf}
We can modify the proof of Theorem 8.9 in \cite{kuo96}.
Use the analyticity in condition (1) to get the expansion
\begin{equation}
F(\x) = \sum_{n=0}^{\infty} \la f_{n}, \x^{\otimes n}\ra.
  \notag
\end{equation}
Then use the Cauchy formula and condition (2) to show that
for any $\x_{1}, \ldots, \x_{n}$ in $[\ce_{p}]_{\a}^{*}$
and any $R>0$:
\begin{equation}
|\la f_{n}, \x_{1}\wh\otimes \cdots \wh\otimes \x_{n}\ra|
\leq {1 \over n!} {KG_{1/\a}(an^{2}R^{2})^{1/2} \over
R^{n}} \,|\x_{1}|_{-p} \cdots |\x_{n}|_{-p}. \notag
\end{equation}
Make a change of variables $r=an^{2}R^{2}$ and use the
inequality $n^{n}\leq n!e^{n}/\sqrt{2\p}$ (see page 357
in \cite{kuo96}) and then take the infimum over $r>0$
to get
\begin{equation}
|\la f_{n}, \x_{1}\wh\otimes \cdots \wh\otimes \x_{n}\ra|^{2}
\leq {K^{2} \over 2\p} a^{n}e^{2n}\left(\inf_{r>0}
{G_{1/\a}(r) \over r^{n}}\right) |\x_{1}|_{-p}^{2} \cdots
|\x_{n}|_{-p}^{2}. \notag
\end{equation}
Now, suppose $q\in [0, p]$ satisfies the condition in
Equation (\ref{eq:2-2}). Then by similar arguments as in
the proofs of Theorems 8.2 and 8.9 in \cite{kuo96} we
can derive
\begin{equation}
|f_{n}|_{q}^{2} \leq {K^{2} \over 2\p} a^{n}e^{2n}
\left(\inf_{r>0} {G_{1/\a}(r) \over r^{n}}\right)
\|i_{p, q}\|_{HS}^{2n}.  \notag
\end{equation}
Therefore
\begin{align}
\|\vf\|_{q, \a}^{2}
  & = \sum_{n=0}^{\infty} n!\a(n) |f_{n}|_{q}^{2} \notag \\
  & \leq {K^{2} \over 2\p} \sum_{n=0}^{\infty}
\left(n!\a(n)\inf_{r>0}{G_{1/\a}(r) \over r^{n}}\right)
\Big(ae^{2}\|i_{p, q}\|_{HS}^{2}\Big)^{n}.  \notag
\end{align}
Note that the series converges because of the condition
in Equation (\ref{eq:2-2}).
\end{pf}

\begin{theorem}  \label{thm:2}
Suppose the sequence $\{\a(n)\}_{n=0}^{\infty}$ satisfies
conditions {\em (A1)}.
\par
\noindent
{\em (I)} Let $\vf\in [\ce]_{\a}$. Then the $S$-transform
$F=S\vf$ of $\vf$ satisfies the conditions:
\begin{itemize}
\item[(1)] For any $\x, \y$ in $\ce_{c}$, the function
$F(z\x+\y)$ is an entire function of $z\in\spc$.
\item[(2)] For any constants $a, p \geq 0$ there exists a
constant $K\geq 0$ such that
\begin{equation} \label{eq:c}
|F(\x)| \leq K G_{1/\a} \big(a|\x|_{-p}^{2}\big)^{1/2},
  \qquad \x\in\ce_{c}.
\end{equation}
\end{itemize}
{\em (II)} Conversely, assume that condition
{\em ($\wt{\text{B}}1$)} holds and let $F$ be a function
on $\ce_{c}$ satisfying the above conditions {\em (1)}
and {\em (2)}. Then there exists a unique $\vf\in
[\ce]_{\a}$ such that $F=S\vf$.
\end{theorem}

We remark that condition (2) is actually equivalent to
the condition: for any constant $p \geq 0$ there exists
a constant $K\geq 0$ such that
\begin{equation}
|F(\x)| \leq K G_{1/\a} \big(|\x|_{-p}^{2}\big)^{1/2},
  \qquad \x\in\ce_{c}.  \notag
\end{equation}
This fact can be easily checked by using the fact that
$|\x|_{-q} \leq \r^{q-p} |\x|_{-p}$ for any $q\geq p$.
Having the constant $a$ in condition (2) is for
convenience to check whether a given function $F$
satisfies the condition.

\begin{pf}
Let $\vf\in [\ce]_{\a}$. Note that $[\ce]_{\a} \subset
[\ce]_{\a}^{*}$ and so $\vf\in [\ce]_{\a}^{*}$. Thus
$F=S\vf$ satisfies condition (1) by Theorem \ref{thm:1}.
To check condition (2), note that
\begin{equation}
F(\x) = S\vf (\x) = \lla :\!e^{\la \cdot, \x\ra}\!:,
\vf \rra,   \qquad \x\in\ce_{c}.  \notag
\end{equation}
Since $\vf\in [\ce]_{\a}$, we have $\|\vf\|_{q, \a}<\infty$
for all $q\geq 0$. Hence
\begin{equation}
|F(\x)| \leq \|\vf\|_{q, \a}\,\|\!:\!e^{\la \cdot,
\x\ra}\!:\!\|_{-q, 1/\a}.  \notag
\end{equation}
But by Equations (\ref{eq:a-1}) and (\ref{eq:a-2})
$\|\!:\!e^{\la \cdot, \x\ra}\!:\!\|_{-q, 1/\a} =
G_{1/\a}\big(|\x|_{-q}^{2}\big)^{1/2}$. Hence
\begin{equation}
|F(\x)| \leq \|\vf\|_{q, \a}\, G_{1/\a}
\big(|\x|_{-q}^{2}\big)^{1/2}. \notag
\end{equation}
For any given $a, p \geq 0$, choose $q\geq p$ such
that $\r^{q-p} \leq \sqrt{a}$. Then
\begin{equation}
|\x|_{-q} \leq \r^{q-p} |\x|_{-p} \leq
\sqrt{a}\, |\x|_{-p}.  \notag
\end{equation}
Therefore, we obtain
\begin{equation}
|F(\x)| \leq \|\vf\|_{q, \a}\,G_{1/\a}\big(|\x|_{-q}^{2}
\big)^{1/2} \leq \|\vf\|_{q, \a}\,G_{1/\a}\big(
a |\x|_{-p}^{2} \big)^{1/2}. \notag
\end{equation}

To prove the converse, assume that condition
($\wt{\text{B}}1$) holds and let $F$ be a function on
$\ce_{c}$ satisfying conditions (1) and (2). Let $q\geq 0$
be any given number. Choose $a, p\geq 0$ such that
\begin{equation}
ae^{2}\|i_{p, q}\|_{HS}^{2}\limsup_{n\to\infty}\left(
n!\a(n) \inf_{r>0}{G_{1/\a}(r) \over r^{n}}\right)^{1/n}
< 1.  \notag
\end{equation}
This inequality can be achieved in two ways: (1) choose any
$a\geq 0$ and then use the fact that $\lim_{p\to\infty}
\|i_{p, q}\|_{HS}=0$, (2) choose any $p\geq 0$ such that
$i_{p, q}$ is a Hilbert-Schmidt operator and then choose a
sufficiently small number $a\geq 0$. (For the first way we
can choose $a=1$ and this is exactly the fact mentioned in
the above remark.) With the chosen $a$ and $p$, use
condition (2) to get a constant $K$ such that the
inequality in Equation (\ref{eq:c}) holds. Then we apply
Lemma \ref{lem:1} to get the inequality in Equation
(\ref{eq:2-3}) so that $\|\vf\|_{q, \a}<\infty$. Hence
$\vf\in [\ce_{q}]_{\a}$ for all $q\geq 0$. Thus $\vf\in
[\ce]_{\a}$ and the converse of the theorem is proved.
\end{pf}

\section{Examples and comments} \label{sec:3}

\noindent
{\bf 1. Four conditions}

In the definition of CKS-space and the characterization
theorems of generalized and test functions we have
assumed several conditions on the sequence
$\{\a(n)\}_{n=0}^{\infty}$:
(A1), (A2), (B1), (B2),
(B3), ($\wt{\text{A}}2$), ($\wt{\text{B}}1$),
($\wt{\text{B}}2$). Recall that we have the following
implications:
\begin{equation}
(\text{A}1) \implies (\wt{\text{A}}2), \quad
(\text{B}2) \implies (\text{B}1), \quad (\wt{\text{B}}2)
\implies (\wt{\text{B}}1), \quad (\text{B}3) \implies
(\wt{\text{B}}2). \notag
\end{equation}
Taking these implications into account we will consider
below the {\bf four conditions}: (A1), (A2), (B2), (B3).

\medskip
\noindent
{\bf 2. Examples}

We give three examples corresponding to the
Hida-Kubo-Takenaka, Kondratiev-Streit, and CKS-spaces.

\begin{example}
For the Hida-Kubo-Takenaka space $(\ce) \subset (L^{2})
\subset (\ce)^{*}$ (see \cite{hkps} \cite{kt80a}
\cite{kt80b} \cite{ob},) the sequence is given by
$\a(n)\equiv 1$. Obviously, this sequence satisfies the
above four conditions. The corresponding exponential
generating functions are
\begin{equation}
G_{\a}(r) = e^{r}, \quad G_{1/\a} (r) = e^{r}. \notag
\end{equation}
Thus the growth conditions in Equations (\ref{eq:d}) and
(\ref{eq:c}) can be stated as
\begin{equation}
|F(\x)| \leq K\,e^{a|\x|_{p}^{2}}, \quad
|F(\x)| \leq K\,e^{a|\x|_{-p}^{2}}.  \notag
\end{equation}
Theorems \ref{thm:1} and \ref{thm:2} are due to
Potthoff-Streit \cite{ps} and Kuo-Potthoff-Streit
\cite{kps}, respectively.
\end{example}

\begin{example} \label{ex:a}
For the Kondratiev-Streit space $(\ce)_{\b} \subset (L^{2})
\subset (\ce)_{\b}^{*}$ (see \cite{ks92} \cite{ks93}
\cite{kuo96},) the sequence is given by $\a(n)=(n!)^{\b}$.
It is easy to check that this sequence satisfies the above
four conditions. The corresponding exponential generating
functions are
\begin{equation}
G_{\a}(r) = \sum_{n=0}^{\infty} {1 \over (n!)^{1-\b}}
  r^{n}, \quad
G_{1/\a}(r) = \sum_{n=0}^{\infty} {1 \over (n!)^{1+\b}}
  r^{n}. \notag
\end{equation}
But from page 358 in \cite{kuo96} and Lemma 7.1 (page 61)
in \cite{cks} we have the inequalities:
\begin{equation} \label{eq:e}
\exp\Big[(1-\b)\,r^{{1\over 1-\b}}\Big] \leq G_{\a}(r)
  \leq 2^{\b}\exp\Big[(1-\b)\, 2^{{\b\over 1-\b}}\,
  r^{{1\over 1-\b}}\Big].
\end{equation}
On the other hand, from page 358 in \cite{kuo96} and the
same argument as in the proof of Lemma 7.1 in \cite{cks}
we can derive the inequalities:
\begin{equation} \label{eq:f}
2^{-\b}\exp\Big[(1+\b)\, 2^{-{\b\over 1+\b}}\,
  r^{{1\over 1+\b}}\Big] \leq G_{1/\a}(r)  \leq
\exp\Big[(1+\b)\,r^{{1\over 1+\b}}\Big].
\end{equation}
The inequalities in Equations (\ref{eq:e}) and (\ref{eq:f})
imply that the growth conditions in Theorems \ref{thm:1} and
\ref{thm:2} are respectively equivalent to the conditions:
\begin{itemize}
\item[$\bullet$] There exist constants $K, a, p \geq 0$
such that
\begin{equation}
|F(\x)| \leq K \exp\Big[a\,|\x|_{p}^{{2\over 1-\b}}
  \Big],  \qquad \x\in\ce_{c}.  \notag
\end{equation}
\item[$\bullet$] For any constants $a, p \geq 0$ there
exists a constant $K\geq 0$ such that
\begin{equation}
|F(\x)| \leq K \exp\Big[a\,|\x|_{-p}^{{2\over 1+\b}}
  \Big],  \qquad \x\in\ce_{c}.  \notag
\end{equation}
\end{itemize}
These two inequalities are the growth conditions used by
Kondratiev and Streit in \cite{ks92} \cite{ks93}
(see also Theorems 8.2 and 8.10 in \cite{kuo96}.)
\end{example}

\begin{example} \label{ex:b}
(Bell's numbers) For each integer $k\geq 2$, consider the
$k$-th iterated exponential function $\exp_{k}(z) =
\exp(\exp(\cdots (\exp(z))))$. This function is entire and
so it has the power series expansion
\begin{equation}
\exp_{k}(z) = \sum_{n=0}^{\infty} {B_{k}(n) \over n!}\,
   z^{n}.  \notag
\end{equation}
The $k$-th order Bell's numbers $\{b_{k}(n)\}_{n=0}^{\infty}$
are defined by
\begin{equation}
b_{k}(n) = {B_{k}(n) \over \exp_{k}(0)},
     \qquad n\geq 0. \notag
\end{equation}
Obviously, this sequence $\{b_{k}(n)\}_{n=0}^{\infty}$
satisfies conditions (A1) and (A2). It has been shown
recently in \cite{akk1} that this sequence satisfies
conditions (B2) and (B3). The corresponding exponential
generating functions are given by
\begin{align}
& G_{b_{k}}(r) = \sum_{n=0}^{\infty} {b_{k}(n) \over n!}
\,r^{n} = {\exp_{k}(r) \over \exp_{k}(0)},
  \label{eq:g}   \\
& G_{1/b_{k}}(r) = \sum_{n=0}^{\infty} {1 \over
n! b_{k}(n)} \,r^{n}.  \label{eq:h}
\end{align}
The function $\exp_{k}(r)$ in Equation (\ref{eq:g}) can be
used to give the growth condition in Theorem \ref{thm:1}
for generalized functions, i.e.,
\begin{itemize}
\item[$\bullet$] There exist constants $K, a, p \geq 0$
such that
\begin{equation}
|F(\x)| \leq K \big(\exp_{k}(a|\x|_{p}^{2})\big)^{1/2},
\qquad  \x\in\ce_{c}.  \notag
\end{equation}
\end{itemize}
On the other hand, we cannot express the sum of the series
in Equation (\ref{eq:h}) as an elementary function. Thus
the growth condition (2) in Theorem \ref{thm:2} is very
hard, if not impossible, to verify. Hence it is desirable
to find similar inequalities as those in Equation
(\ref{eq:f}) for the sequence $\{b_{k}(n)\}_{n=0}^{\infty}$.
\end{example}

\noindent
{\bf 3. General growth order}

Being motivated by Examples \ref{ex:a} and \ref{ex:b}, we
consider the question: What are the possible functions $U$
and $u$ such that the growth conditions in Theorems
\ref{thm:1} and \ref{thm:2} can be respectively stated
as follows?
\begin{itemize}
\item[$\bullet$] There exist constants $K, a, p \geq 0$
such that
\begin{equation}
|F(\x)| \leq K U\big(a\,|\x|_{p}^{2}\big)^{1/2},
   \qquad \x\in\ce_{c}.  \notag
\end{equation}
\item[$\bullet$] For any constants $a, p \geq 0$ there
exists a constant $K\geq 0$ such that
\begin{equation}
|F(\x)| \leq K u\big(a\,|\x|_{-p}^{2}\big)^{1/2},
  \qquad \x\in\ce_{c}.  \notag
\end{equation}
\end{itemize}

The answer to this question will be given in our forthcoming
papers \cite{akk2} \cite{akk3}. In particular, when $\a(n)
=b_{k}(n)$, condition (2) in Theorem \ref{thm:2} can be
replaced by the following growth condition:
\begin{itemize}
\item[$\bullet$] For any constants $a, p \geq 0$ there
exists a constant $K\geq 0$ such that
\begin{equation}
|F(\x)| \leq K \exp\left[\sqrt{a|\x|_{-p}^{2}\,
\log_{k-1} \big(a|\x|_{-p}^{2}\big)}\,\right],
  \qquad \x\in\ce_{c},  \notag
\end{equation}
\end{itemize}
where $\log_{j}, \,j\geq 1,$ is the function defined by
\begin{equation}
\log_{1}(x)=\log(\max\{x, e\}), \quad \log_{j}(x) =
\log_{1}(\log_{j-1}(x)), \quad j\geq 2.  \notag
\end{equation}

\end{document}